\magnification=1200
\baselineskip=14pt

\def\qed{{$\vrule height4pt depth0pt width4pt$}}
\def\ss{{\smallskip}}
\def\ms{{\medskip}}
\def\bs{{\bigskip}}
\def\ni{{\noindent}}
\def\S{{\T^{\Z^d}}}

\def\s{{\sigma}}
\def\a{{\alpha}}
\def\t{{\tau}}

\def\g{{\gamma}}
\def\d{{\delta}}
\def\b{{\beta}}
\def\Z{{\bf Z}}

\def\<{{\langle}}
\def\>{{\rangle}}

\def\D{{\Delta}}
\def\T{{\bf T}}
\def\R{{\bf R}}
\def\<{{\langle}}
\def\>{{\rangle}}
\def\w{{\wedge}}
\def\m1{{\quad ({\rm mod\ 1})}}
\def\L{{\Lambda}}
\def\m{{\bf m}}
\def\Rd{{\cal R}_d}
\def\Sd{{\bf S}^d}
\input epsf
\epsfverbosetrue

\centerline{\bf MAHLER MEASURE, LINKS AND HOMOLOGY GROWTH}
\smallskip
\centerline{Daniel S. Silver and Susan G. Williams} \ss

\ms
\footnote{} {Both authors partially supported by NSF grant DMS-9704399.}
\footnote{}{2000 {\it Mathematics Subject Classification.}  Primary 57M25; secondary 37B10.}

\noindent {\narrower\narrower\smallskip\noindent  {\bf ABSTRACT.} Let $l$ be an oriented
link of $d$ components with nonzero Alexander polynomial $\D(u_1, \ldots, u_d)$. Let $\L$ be 
a finite-index subgroup of $H_1(S^3-l)\cong \Z^d$, and let $M_\L$ be the corresponding
abelian cover of $S^3$ branched  along $l$.  The growth
rate of the order of the torsion subgroup of
$H_1(M_\L)$, as a suitable measure of $\L$ approaches infinity, is equal to the Mahler measure of $\Delta$.
\smallskip}
\ms

\ni {\bf 1. Introduction.} Associated to any knot $k\subset S^3$ is a sequence of Alexander polynomials
$\D_i,\ i\ge 1$, in a single variable such that $\D_{i+1}$ divides $\D_i$. Likewise, for any oriented
link of $d$ components there is a sequence of Alexander polynomials in $d$ variables. Following the
usual custom, we refer to the first Alexander polynomial of a knot or a link as {\sl the}
Alexander polynomial, and we denote it simply by $\D$. 

In [{\bf Go}] C.~McA. Gordon examined the homology groups of $r$-fold
cyclic covers $M_r$ of $S^3$ branched over a knot $k$. He proved that when each zero of 
the Alexander polynomial $\D$ of $k$ has modulus one (and hence is a root of unity), the
finite values of $|H_1(M_r)|$ are periodic in $r$.  Gordon conjectured that when some zero of $\D$ is not a
root of unity, the finite values of $|H_1(M_r)|$ grow exponentially.

More than fifteen years later two independent proofs of Gordon's conjecture, one by R. Riley [{\bf Ri}] and 
another by F. Gonzal\'ez-Acu\~na and H. Short [{\bf GoSh}], appeared. Both employed the
Gel'fond-Baker theory of linear forms in the logarithms of algebraic integers [{\bf Ba}],
[{\bf Ge}]. 

We extend the above results for knots, replacing the term ``finite values of $|H_1(M_r)|$'' with 
``order of the torsion subgroup of $H_1(M_r)$,'' while at the same time proving a general result for 
links in $S^3$. Our proof, which is motivated by [{\bf SiWi2}], identifies the torsion
subgroup of  the homology of a finite abelian branched cover with the connected components of 
periodic points in an associated algebraic dynamical system. 
Theorem 21.1 of [{\bf Sc}], an enhanced version of a theorem of D. Lind, K. Schmidt and T. Ward [{\bf
LiScWa}], then completes our argument. 

Recognizing that relatively few topologists are familiar with algebraic dynamical systems, we have
endeavored to make this paper self-contained. The reader who desires to know more about such dynamical
systems is encouraged to consult the extraordinary monograph [{\bf Sc}].

We thank Jonathan Hillman, Douglas Lind and Klaus Schmidt for valuable discussions. 
Also, we are grateful to the University of Washington for its hospitality
when part of this work was completed. 

\bs

\ni {\bf 2. Statement of results.} Let $l= l_1 \cup \cdots \cup l_d$ be an oriented link of 
$d$ components with exterior $E = S^3 - {\rm int}N(l)$, where $N(l)$ is a regular neighborhood of $l$. The
meridianal generators of the link group $G_l= \pi_1(S^3 -l)$ represent  distinguished generators 
$u_1, \ldots, u_d$ for the abelianization $G_l/G_l' \cong \Z^d$. We identify these generators with the 
standard basis of $\Z^d$. 

Given  a finite-index subgroup 
$\L \subset \Z^d$ there exists a covering $E_\L$ of the link exterior corresponding to the epimorphism
$G_l \to \Z^d/\L$, the abelianization map composed with the canonical quotient map. By attaching solid
tori to
$E_\L$ so that meridians of the tori  cover meridians of $l$ while the collection of cores map to the link,
we obtain a cover $M_\L$ of
$S^3$ branched over $l$. 

Let $\Rd$ denote the ring $\Z[u_1^{\pm 1}, \ldots, u_d^{\pm 1}]\cong \Z[\Z^d]$ of Laurent polynomials with
integer coefficients. The {\it Mahler measure} of a nonzero polynomial $f \in \Rd$ is defined by
$${\bf M}(f) = \exp (\int_{\Sd}\log |f({\bf s})|\ d{\bf s}),$$
where $d{\bf s}$ indicates integration with respect to normalized Haar measure, and $\Sd$ is the 
multiplicative subgroup of $d$-dimensional complex space ${\bf C}^d$ consisting of all vectors
$(s_1, \ldots, s_d)$ with $|s_1|= \cdots =|s_d| = 1$. Clearly, Mahler measure is multiplicative, and
the measure of any unit is $1$. It is known that $M(f)=1$ if and only if $f$ is equal up to a unit factor
to the product of cyclotomic polynomials in a single variable evaluated at monomials 
(see [{\bf Sc}, Lemma 19.1]). 

The quantity ${\bf M}(f)$, which is the geometric mean of $|f|$ over the $d$-torus $\Sd$,
was introduced by K.~Mahler in [{\bf Ma1}] and [{\bf Ma2}]. It is a consequence of Jensen's formula [{\bf
Al}, p. 208] that if $f$ is a nonzero polynomial $c_nu^n + \cdots +c_1u +c_0\ (c_n\ne 0)$ in one variable,
then
$${\bf M}(f) = |c_n|\cdot \prod_{j=1}^n {\rm max}(|r_j|,1),$$
where $r_1,\ldots, r_n$ are the zeros of $f$. A short proof can be found in either [{\bf EvWa}] or [{\bf Sc}]. 

For any finite-index subgroup $\L$ we 
let
$$\<\L\> = {\rm min}\{|v|\ :\ v \in \L - 0\},$$
where $|\cdot|$ denotes the Euclidean metric. 

Since $M_\L$ is a compact manifold, the homology group $H_1(M_\L)$ is finitely generated. (All homology
groups in this paper have integer coefficients.) We decompose $H_1(M_\L)$ as the direct
sum of a free abelian group of some rank $\b_\L$ and a torsion subgroup $TH_1(M_\L)$. We denote
the order of $TH_1(M_\L)$ by $b_\L$. \bs

\ni {\bf Theorem 2.1.} Let $l=l_1\cup \ldots \cup l_d$ be an oriented link of $d$ components
having nonzero Alexander polynomial $\D= \D(u_1,\ldots, u_d)$. Then
$$\limsup_{\<\L\> \to \infty}\ {1 \over {|\Z^d/\L|}} \log b_\L = \log {\bf M}(\D).$$
When $d=1$ the $\limsup$ can be replaced by an ordinary limit. \bs

As a consequence of the proof of Theorem 2.1 we obtain a new proof of a theorem
of Gordon. Recall that for any knot, $\D_i$ denotes the $i$th Alexander polynomial. \bs

\ni {\bf Corollary 2.2.} [{\bf Go}] Let $k$ be a knot in $S^3$. If $\D_1/\D_2$ divides
$t^N-1$ for some $N$, then
$H_1(M_r)\cong H_1(M_{r+N})$ for all $r$. 

\bs
The Mahler measure of $u_1^2-u_1+1$, the Alexander polynomial of the 
trefoil knot $3_1$, is $1$ since both zeros of the polynomial have unit modulus. On the other
hand, $u_1^2-3u_1+1$, the Alexander polynomial of the figure eight knot $4_1$, has zeros
$(1 \pm \sqrt 5)/2$ and hence it has Mahler measure $(1+\sqrt 5)/2 \approx 1.618$. 

Scanning the table of $2$-component links in [{\bf Ro}] we find that $6_2^2$ is the first
link with nonzero Alexander polynomial having Mahler measure greater than $1$.  The
polynomial is $u_1 +u_2 -1 + u_1^{-1} + u_2^{-1}$, which has Mahler measure approximately
equal to $1.285$.  

The next link in the table, $6_2^3$, has Alexander polynomial $2-u_1-u_2+2u_1u_2$, which
can be rewritten as $(2-u_1)+u_1u_2(2u_1-1)$. Using Lemma 19.8 of [{\bf Sc}] and an easy 
change of basis (replacing $u_1u_2$ with a new variable $u_2'$) we see that 
the Mahler measure of this Alexander polynomial is precisely $2$.

No polynomial with integer coefficients is known that has Mahler measure greater than 1 but less than that
of $\D(x) = x^{10}+x^9-x^7-x^6-x^5-x^4-x^3+x+1$, which has
Mahler measure approximately equal to $1.176$. (Only one of the
nine zeros of
$\D$ lies outside the unit circle.) Deciding whether or
not such a numerical gap truly exists is known as Lehmer's
problem, and it remains a vexing open question. (see [{\bf Le}],
[{\bf EvWa}]). It is a provocative fact that $\D$ is the
Alexander polynomial of a knot. In fact, there are infinitely
many, including infinitely many with complements that fiber over
the circle.

For more calculations of Mahler measures of Alexander polynomials of links and
further discussion of Lehmer's question see [{\bf SiWi4}].

\bs

\ni {\bf 3. Z$^d$-shifts and link colorings.} 
A $\Z^d$-{\it action by automorphisms} on a topological group
$X$ is a homomorphism $\s: \m \mapsto \s_\m$ from $\Z^d$ to ${\rm Aut}(X)$.  Two
$\Z^d$-actions $\s$ and $\s'$ on $X$ and $X'$, respectively, are {\it
algebraically conjugate} if there exists a continuous group isomorphism
$\phi: X \to X'$ such that $\phi\circ \s_\m = \s'_\m\circ \phi$, for every $\m \in \Z^d$.

$\Rd$-modules are an important source of $\Z^d$-shifts via Pontryagin duality. Let $\T$ denote
the additive circle group $\R/\Z$. For any $\Rd$-module
$L$, the Pontryagin dual $L^\wedge={\rm Hom}(L, \T)$ is a group under pointwise
addition. Here $L$ is given the discrete topology, and
$L^\wedge$ is endowed with the  compact-open topology; $L^\wedge$ is a compact abelian group.
For $\m \in
\Z^d$, scalar multiplication $a \mapsto \m a$ in $L$ induces $\s_\m \in {\rm Aut}(L^\wedge)$ via its
adjoint action. In
this way we have a $\Z^d$-action on $L^\wedge$. From a purely algebraic point of view, 
$L^\wedge$ is a $\Rd$-module. 
In the case that $L$ is a free $\Rd$-module of rank $N$, we obtain the compact group $L^\wedge = (\T^N)^{\Z^d}$.
The automorphism $\s_\m$ is the shift map given by $\s_\m(\a_{\bf n}) = (\a_{{\bf n}+{\bf m}})$ for 
$\a= ( \a_{\bf
n}) \in (\T^N)^{\Z^d}$. The automorphisms
$\s_{u_1}, \ldots, \s_{u_d}$ will be denoted by $\s_1, \ldots, \s_d$ for notational ease. 

Given a $\Z^d$-action $\s$ on $X$, we say that a point $x \in X$ is {\it periodic}
under $\s$ if its orbit $\{\s_\m x\ |\ \m \in \Z^d\}$ is finite. We will be particularly interested
in periodic point sets
$${\rm Fix}_\L(\s) = \{x \in X\ |\ \s_\m x=x\ \forall \m\in \L \},$$
where $\L$ is a subgroup of finite index in $\Z^d$. For algebraically conjugate actions such sets
clearly correspond under the isomorphism $\phi$. \bs

\ni {\bf Definition 3.1.} Assume that $D$ is a diagram of an oriented link $l=l_1 \cup \cdots \cup l_d$ of
$d$ components. A $\S$-{\it coloring} of $D$ is an assignment of elements ({\it colors}) $\a, \b, \ldots \in
\S$ to the  arcs of $D$ such that the condition 
$$\a + \s_t \b = \g + \s_{t'}\a \eqno(3.1)$$
is satisfied at any crossing. Here $\a$ corresponds to an overcrossing arc of the $t$th component
of
$l$, while $\b$ and $\g$ correspond to undercrossing arcs of the $t'$th component. We encounter
$\b$ as we travel in the preferred direction along the arc labeled by $\a$, turning left at the crossing
(see Figure 1). The terminology is motivated by the concept of Fox coloring for knots [{\bf Fo}], which was
generalized in [{\bf SiWi1}], [{\bf SiWi2}].\bs

\epsfxsize=.8truein
\centerline{\epsfbox{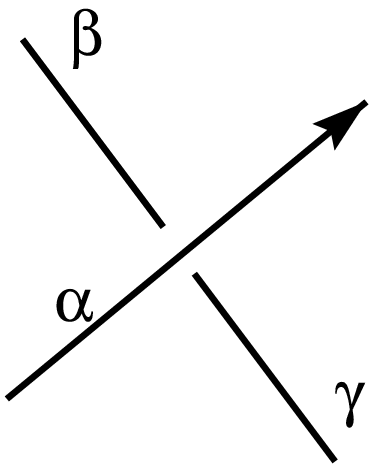}}
\bs
\centerline{{\bf Figure 1:} $\S$-Coloring rule}
\bs
If $D$ consists of $N$ arcs, then the set ${\rm Col}_{\T, \Z^d}(D)$ of all $\S$-colorings of $D$ is
a closed subgroup of $[\T^{\Z^d}]^N \cong [\T^N]^{\Z^d}$ that is 
invariant under $\s_\m$ for each $\m \in \Z^d$. It will follow from observations in Section 4 that if $D'$ is a another diagram for $l$, then the
$\Z^d$-action on
${\rm Col}_{\T,
\Z^d}(D)$ is algebraically conjugate to the action on ${\rm Col}_{\T, \Z^d}(D')$. In anticipation we
make the following definition. \bs

\ni {\bf Definition 3.2.} (Cf. [{\bf SiWi1}]) Let $l$ be a $d$-component oriented link with diagram $D$.
The {\it color $\Z^d$-shift} ${\rm Col}_{\T,\Z^d}(l)$ is the compact abelian group 
${\rm Col}_{\T, \Z^d}(D)$ together with the $\Z^d$-action $\s$. \bs

\ni {\bf 4. Alexander module and periodic points.} This section contains the proofs of our main results.
An example that illustrates the main ideas is given at the end.  

A diagram $D$ for $l$ yields a finite Wirtinger
presentation 
$\<x_1, x_2 \ldots, x_N \mid r_1, \ldots, r_N\>$ for $G_l = \pi_1(S^3-l)$. 
Let $P$ be the canonical $2$-complex with $\pi_1 P \cong G_l$, constructed with a single vertex $v$,
directed edges labeled $x_1, \ldots, x_N$, and oriented $2$-cells
$c_1, \ldots, c_N$ with each boundary $\partial c_i$ attached to $1$-cells according to $r_i$
(see Chapter 11 of [{\bf Li}]). 

Let $\tilde P$ be the maximal abelian cover of $P$; that is, the cover corresponding to the abelianization
map $G_l \to G_l/G_l' \cong \Z^d$. As usual each cell $v, x_i, c_j$ of $P$ lifts to a family $\m \tilde
v,\
\m \tilde x_i,\ \m \tilde c_j$ of oriented cells indexed by $\Z^d$. By standard construction the
chain complex $0 \to \tilde C_2\ {\buildrel \partial_2 \over \to}\ \tilde C_1 {\buildrel \partial_1 \over
\to}\ \tilde C_0
\to 0$ admits a quotient $0 \to \tilde C_2\ {\buildrel \partial_2
\over \to}\ \tilde C_1 \to 0$ that determines the relative homology group $H_1(\tilde P, \tilde P^0)$, 
where $\tilde P^0$ is the $0$-skeleton of $\tilde P$. This $\Rd$-module is the Alexander module of 
the link, denoted here by $A$. 

By the universal coefficient theorem [{\bf Sp}, p. 243] the cohomology group $H^1(\tilde P, \tilde P^0;
\T)$ is naturally isomorphic to the dual group of the Alexander module. It is a closed subgroup of
${\rm Hom}(\tilde C_1, \T) = [\T^N]^{\Z^d}$, the kernel of the coboundary operator $${\rm Hom}(\tilde C_1,\T)
\buildrel {\rm Hom}(\partial_2,1)
\over \longrightarrow {\rm Hom}(\tilde C_2, \T),$$ and hence it inherits a $\Z^d$-action from ${\rm
Hom}(\tilde C_1,\T)$ (see Section 3). 

We observe that
$\partial_2 \tilde c_1, \ldots, \partial_2 \tilde c_N$ closely resemble the relations of the coloring rule
(3.1). The Wirtinger relator at the crossing in Figure 1 has the form $x_ix_{j_1}x_i^{-1}x_{j_2}^{-1}$, and
the lifted loop in the cover that begins at $\tilde v$ determines the $1$-cycle
$\tilde x_i + u_t \tilde x_{j_1} -u_{t'} \tilde x_i - \tilde x_{j_2}$, which induces the homology relation
$\tilde x_i + u_t \tilde x_{j_1} = \tilde x_{j_2} + u_{t'} \tilde x_i$. Lifts that 
begin at other points of the cover are simply translates by elements of $\Z^d$. We can regard
the assignment of $\a=(\a_\m)\in \T^{\Z^d}$ to an arc $x_i$ as an assignment of $\a_\m \in \T$ to the $1$-chain 
$\m \tilde x_i$. Then clearly ${\rm Col}_{\T, \Z^d}(D)$ and 
$H^1(\tilde P, \tilde P^0;\T)$ are described by identical subsets of $[\T^N]^{\Z^d}$. 
Since the Alexander module is a link invariant, it follows that the algebraic conjugacy class of
${\rm Col}_{\T, \Z^d}(D)$ is independent of the diagram for $l$.

Consider a finite-index subgroup $\L$ of $\Z^d$. The unbranched cover $E_\L$ has the same homology
as the quotient complex $\tilde P/\L$. A $2$-complex $Q$ with the same first homology group as
the branched cover $M_\L$ is obtained from $\tilde P/\L$ by attaching additional $2$-cells
as follows. Each Wirtinger generator
$x_i,\ 1\le i \le N$, maps to some $u_{t(i)}$ under abelianization. Assume that $u_{t(i)}$
represents an element of  order $n(i)$ in $\Z^d/\L$. Then $\tilde x_i +
u_{t(i)} \tilde x_i + \cdots + u_{t(i)}^{n(i)-1}\tilde x_i$ is a  $1$-cycle in $\tilde P/\L$, and we attach
a $2$-cell along it and each of its translates. The additional cells added to the complex in this way correspond
to the meridinal disks of tori  that we attach to $E_\L$ when constructing $M_\L$. Consequently, $Q$ has the same
fundamental group and  hence the same first homology group as $M_\L$.

Elements of $H^1(Q,Q^0;\T)$ correspond to 
$\Z^d$-colorings in ${\rm Fix}_\L(\s)$ such that if $\a$ is assigned to the $i$th arc of $D$, then 
$$\a+\s_{t(i)}\a+\cdots+\s_{t(i)}^{n(i)-1}\a=0.\eqno (4.1)$$ 
Such $\Z^d$-colorings comprise a subgroup 
${\rm SFix}_\L(\s)\le{\rm Fix}_\L(\s)$ of {\it special periodic points}. 

The universal coefficient theorem implies that $H^1(Q;\T)$ is isomorphic to the dual group of $H_1(Q) \cong
H_1(M_\L)$.  By decomposing $H_1(M_\L)$ as $TH_1(M_\L) \oplus \Z^{\b_\L}$ (see Section 2) and recalling
that $A^\w \cong A$ for any finite abelian group $A$, we have
$$H^1(Q;\T) \cong [TH_1(M_\L) \oplus \Z^{\b_\L}]^\w \cong TH_1(M_\L)\oplus \T^{\b_\L}.\eqno (4.2)$$ 

Consider now the portion of the cohomology long exact sequence:
$$H^0(Q^0; \T) {\buildrel \d \over \to} H^1(Q,Q^0;\T) \to H^1(Q;\T)
\to 0.\eqno(4.3)$$

\ni {\bf Lemma 4.1.} The image of $\d$ is a direct summand of $H^1(Q,Q^0;\T)$ isomorphic to $\T^r$, 
where $r = |\Z^d/\L|-1$. \bs

\ni {\bf Proof.} The $0$-skeleton $Q^0$ consists of vertices indexed by elements of $\Z^d/\L$. Elements of
$H^0(Q^0;\T)$ can be regarded as functions $f: \Z^d/\L \to \T$. The image $\d(f)$ is an edge-labeling of $Q$, assigning $f(\m')-f(\m)$ to an edge
from $\m \tilde v$ to $\m' \tilde v$. Select a maximal tree $T$ in the $1$-skeleton of $Q$. It is clear that
$\d(f)$ is uniquely determined by its values on $T$, and such values can be
prescribed arbitrarily.  Since $T$ has $r$ edges, the image of $\d$ is isomorphic to $\T^r$.

The map $H^1(Q,Q^0;\T) {\buildrel \epsilon \over \to} \T^r$ given by restricting any cocycle to the 
edges of the maximal tree $T$ is an epimorphism, by what we have said above. We construct a
right inverse $\eta$ for $\epsilon$ as follows. Given a function $g: T \to \T$, choose an element $f \in
H^0(Q^0;\T)$ such that $\d(f)$ agrees with $g$ on $T$. Define $\eta(g) = \d(f)$. Hence the image of $\d$ is
a direct summand of $H^1(Q,Q^0;\T)$. \qed \bs

\ni {\bf Corollary 4.2.} Let ${\rm SFix}_\L^0(\s)$ be
the connected component of the identity in 
${\rm SFix}_\L(\s)$. Then ${\rm SFix}_\L(\s)/{\rm SFix}^0_\L(\s) \cong TH_1(M_\L)$. \bs

\ni {\bf Proof.} By equation (4.2) we have $TH_1(M_\L)\oplus \T^{\b_\L} \cong H^1(Q;\T)$, and by the
long exact sequence (4.3) the latter module is isomorphic to $H^1(Q,Q^0;\T)/{\rm im}(\d)$. Recall that
$H^1(Q,Q^0;\T)$ is isomorphic to
${\rm SFix}_\L(\s)$. By Lemma 4.1 the image of $\d$ is connected and hence contained in ${\rm
SFix}_\L^0(\s)$. Thus ${\rm SFix}_\L(\s)/{\rm SFix}_\L^0(\s) \cong TH_1(M_\L)$.  \qed \bs

\ni {\bf Corollary 4.3.} The first Betti number $\b_\L$ of $M_\L$ is equal to 
${\rm dim\ SFix}_\L(\s)-|\Z^d/\L|+1.$ \bs

\ni {\bf Proof.} By equation (4.2) the Betti number $\b_\L$ is equal to the dimension of $H^1(Q;\T)$,
and by the long exact sequence (4.3) the latter is ${\rm dim\ SFix}_\L(\s) - {\rm dim\ im}(\d)$. Lemma
4.1 completes the argument. \qed \bs

The quantity ${\rm dim\ SFix}_\L(\s)$ can be generally computed as the nullity of a certain matrix. 
The computation is similar to that which results from the formula of M. Sakuma [{\bf Sa}, Theorem
1.1(2)]. However, the dynamical systems perspective here is new.  \bs

\ni {\bf Proof of Theorem 2.1.} It follows from Corollary 4.2 that $b_\L = |TH_1(M_\L)|$ is equal to the number
of connected components of ${\rm SFix}_\L(\s)$. Now we apply techniques of symbolic dynamics to count
the number of components and determine their exponential growth rate.

We denote the connected component of the identity in ${\rm Fix}_\L(\t)$ by the symbol ${\rm Fix}^0_\L(\t)$.
For any
$\Z^d$-action $\t$ associated to the dual group of a
$\Rd$-module, Theorem 21.1 of [{\bf Sc}] implies that  the exponential growth
rate of $|{\rm Fix}_\L(\t)/{\rm Fix}_\L^0(\t)|$ as $\<\L\>$ approaches infinity is equal to the topological
entropy of
$\t$, {\sl provided that the topological entropy of $\t$ is finite}. The topological entropy
of the $\Z^d$-action
$\s$ above is always infinite. However, we will show that there is a related $\Z^d$-shift $\s'$ such
that (1) $\s'$ has finite topological entropy equal to $\log {\bf M}(\D)$; and (2)
$|{\rm Fix}_\L(\s')/{\rm Fix}_\L^0(\s')| = |{\rm SFix}_\L(\s)/{\rm SFix}_\L^0(\s)|$. The main conclusion
of Theorem 2.1 follows from these assertions. 

A $\T^{\Z^d}$-coloring of $D$ that assigns $0 \in \S$ to some arc, say the arc corresponding to Wirtinger 
generator $x_1$, will be called a {\it based}
$\T^{\Z^d}$-coloring. The collection of based $\T^{\Z^d}$-colorings is a closed shift-invariant subgroup
of ${\rm Col}_{\T, \Z^d}(D)$, independent of the choice of arc.  It
is the dual of $B = A/\<\tilde x_1\>$, the quotient of the Alexander module  by the
submodule generated by $\tilde x_1$, which one might call the {\it based Alexander module} of the link.  We
denote the
$\Z^d$-action on $B$ by
$\s^B$.

It is clear from the discussion above that for any knot the based Alexander module is isomorphic to the first
homology group of the infinite cyclic cover of the knot exterior. Since all generators of a Wirtinger
presentation for a knot group are conjugate, it follows that periodic points of $\s^B$ are always special
periodic points. (See the paragraph preceding Lemma 2.7 [{\bf SiWi3}] for details.) Hence in this case
${\rm Fix}((\s^B)^r)$ is isomorphic to the dual of
$H_1(M_r)$, for any $r$.  

In the general case we define
special periodic points of $\s^B$ just as we defined them for $\s$. 
We claim that 
$${\rm SFix}_\L(\s^B)/{\rm SFix}^0_\L(\s^B) \cong {\rm SFix}_\L(\s)/{\rm SFix}^0_\L(\s),$$
and hence the number of connected components of
${\rm SFix}_\L(\s^B)$ is equal to that of ${\rm SFix}_\L(\s)$. One way to see this is by
constructing a maximal tree $T$ for the $1$-skeleton of $Q$, selecting first a maximal number of edges of
the form $\m
\tilde x_1$. Recall that ${\rm SFix}_\L(\s)$ can be identified with $H^1(Q, Q^0;\T)$, and by the proof of
Lemma 4.1 this group has a direct summand $\T^r$, the image of the map $\eta$ defined above. 
We pass from ${\rm SFix}_\L(\s)$ to ${\rm SFix}_\L(\s)/{\rm SFix}^0_\L(\s)$ in two stages.
First we crush $\eta (\T^r)$; the quotient group is isomorphic to the subgroup of $H^1(Q,
Q^0;\T)$ consisting of cocyles that vanish on the edges of the maximal tree $T$.  However, in view
of equation (4.1), any cocycle that vanishes on those edges vanishes on all of the edges
$\m\tilde x_1$. Thus ${\rm SFix}_\L(\s)/\eta (T^r) \cong {\rm SFix}_\L(\s^B)$,  the subgroup of $H^1(Q, Q^0; \T)$
consisting of  based $\T^{\Z^d}$-colorings. The connected component of the identity in this group is
also a torus.  Crushing it we find that  ${\rm SFix}_\L(\s)/{\rm SFix}_\L^0(\s)\cong {\rm SFix}_\L(\s^B)/{\rm
SFix}^0_\L(\s^B)$.

The based Alexander module $B$ has an $(N-1)\times (N-1)$ presentation matrix $R$, which can be obtained
from the matrix for the Alexander module $A$ by deleting the first row and column. Then
$\D(u_1,\ldots, u_d)=(u_1 - 1)\ {\rm det}\ R$ (see [{\bf Li}], pp. 119 -- 120). Since the Mahler
measure of $u_1 - 1$ is equal to $1$, the determinant of $R$ has the same Mahler measure as
$\D$. By [{\bf LiScWa}, p. 611] the topological entropy of $\s^B$ is equal to the log of 
the Mahler measure of $\D$. 

The $\Z^d$-action $\s'$ that we need is a modification of
$\s^B$. Consider based $\Z^d$-colorings of $D$, but replace any color $\b$ by a pair $(\b, \zeta)$ of colors,
and require in addition to the basic coloring rule (3.1) the condition: $$\s_t\zeta = \zeta +
\b,\eqno(4.4)$$ where
$t$ is the index of the  component of $l$ containing the arc colored by $(\b,\zeta)$. Denote the associated
$\Z^d$-action by $\s'$.

The $\Z^d$-action $\s'$ is on the dual group of a module $B'$ that we obtain from a presentation for
$B$ by adding new generators $z_2, \ldots, z_N$ and relations $u_{t(i)}z_i = z_i + x_i$ ($2\le i\le N$),
where $t(i)$ is the index of the component of $l$ containing the arc $x_i$. The determinant
of the new relation matrix $R'$ is $\Delta(u_1, \ldots, u_d)$ times a product of
polynomials of the form $u_t -1$. As before, since the Mahler measure of each $u_t-1$ is equal to $1$, we
have ${\bf M}({\rm det}(R')) = {\bf M}(\Delta)$, and by [{\bf LiScWa}, p. 611] the topological
entropy of $\s'$ is equal to this value. Hence assertion (1) above holds. 

By a straightforward recursion argument we find that ${\rm Fix}_\L(\s') \cong 
{\rm SFix}_\L(\s^B)\oplus \T^s$. Here $s$ is the number of second coordinates $\zeta$
that can be freely assigned: Assume that $u_t$ represents an element of order $n$ in $\Z^d/\L$.
Condition (4.4) implies
$$\zeta_{\m+u_t}= \zeta_\m +\b_\m$$
$$\zeta_{\m+2u_t} = \zeta_\m+ \b_\m +\b_{\m+u_t}$$
\line{\hfil \vdots \hfil}
$$\zeta_{\m+nu_t} = \zeta_\m+\b_\m + \b_{\m+u_t}+\cdots +\b_{\m+(n-1)u_t}.$$
Clearly, $\zeta_{\m+nu_t} = \zeta_\m$ if and only if $\b_{\m+u_t}+\cdots
+\b_{\m+(n-1)u_t}=0.$ Moreover, the coordinates $\zeta_{u_t},\ldots,
\zeta_{(n-1)u_t}$ are uniquely determined from $\zeta_{\bf 0}$ and coordinates of $\b$.
When $t'$ is different from $t$, condition (4.4) imposes no new requirement; in such a case the
coordinates $\zeta_{\bf 0}, \zeta_{u_{i'}}, \ldots, \zeta_{(n'-1)u_{t'}}$ can be chosen
arbitrarily, where $u_{t'}$ represents an element of order $n'$ in $\Z^d/\L$.  Assertion
(2) is immediate, and the proof of the theorem is complete.  
\qed \bs

\ni {\bf Example 4.4.} We will illustrate the ideas and terminology used in the proof of Theorem 2.1
with an example. 

The diagram $D$ for the link $l=5_1^2$, shown in Figure 2, yields a Wirtinger 
presentation for $G_l$:
$$\<x_1,x_2,x_3,x_4,x_5\mid x_1x_3=x_5x_1, x_3x_2=x_1x_3, x_5x_4=x_3x_5, x_4x_2=x_1x_4, x_2x_4=x_5x_2\>.$$
\epsfxsize=2truein
\centerline{\epsfbox{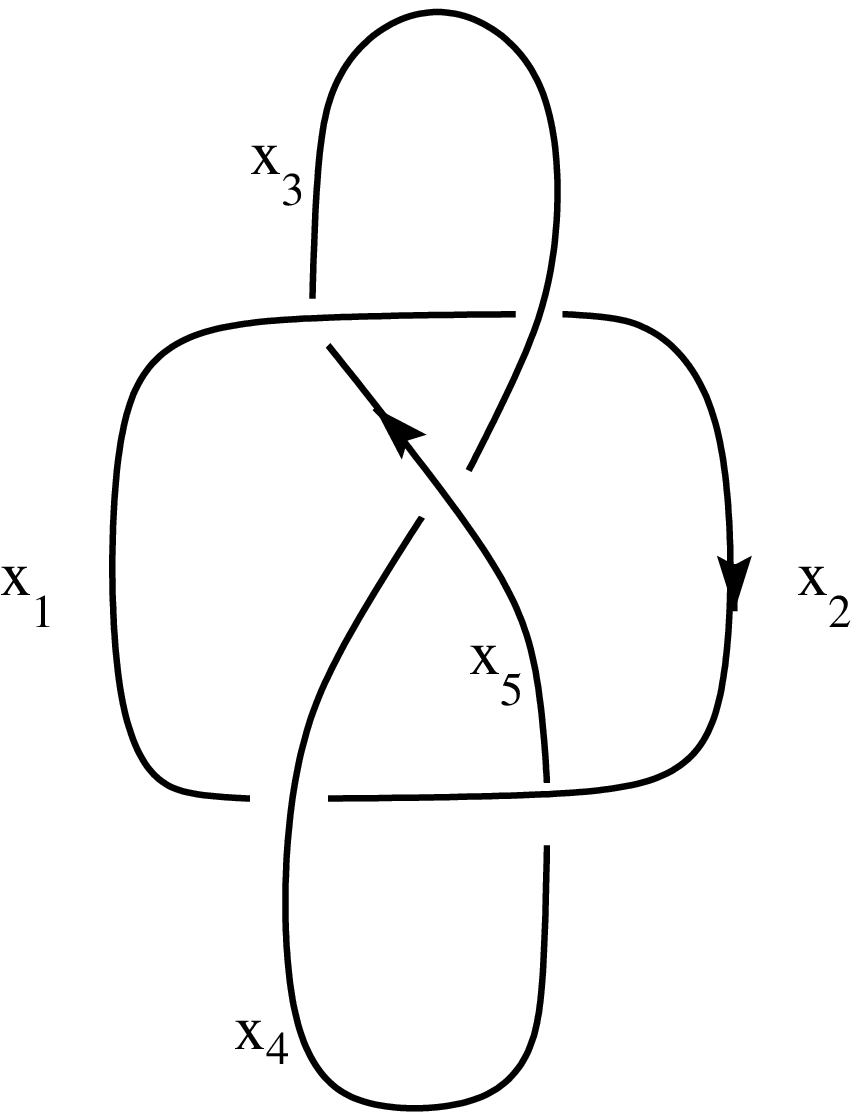}}
\bs
\centerline{{\bf Figure 2:} Diagram for $5_1^2$}
\bs

We assume that under abelianization $x_1$ and $x_2$ map to $u_1$ while the remaining generators are sent to
$u_2$. A portion of the maximal abelian cover $\tilde P$ is shown in Figure 3. The
$2$-cells are not shown. \bs

\epsfxsize=3.2truein
\centerline{\epsfbox{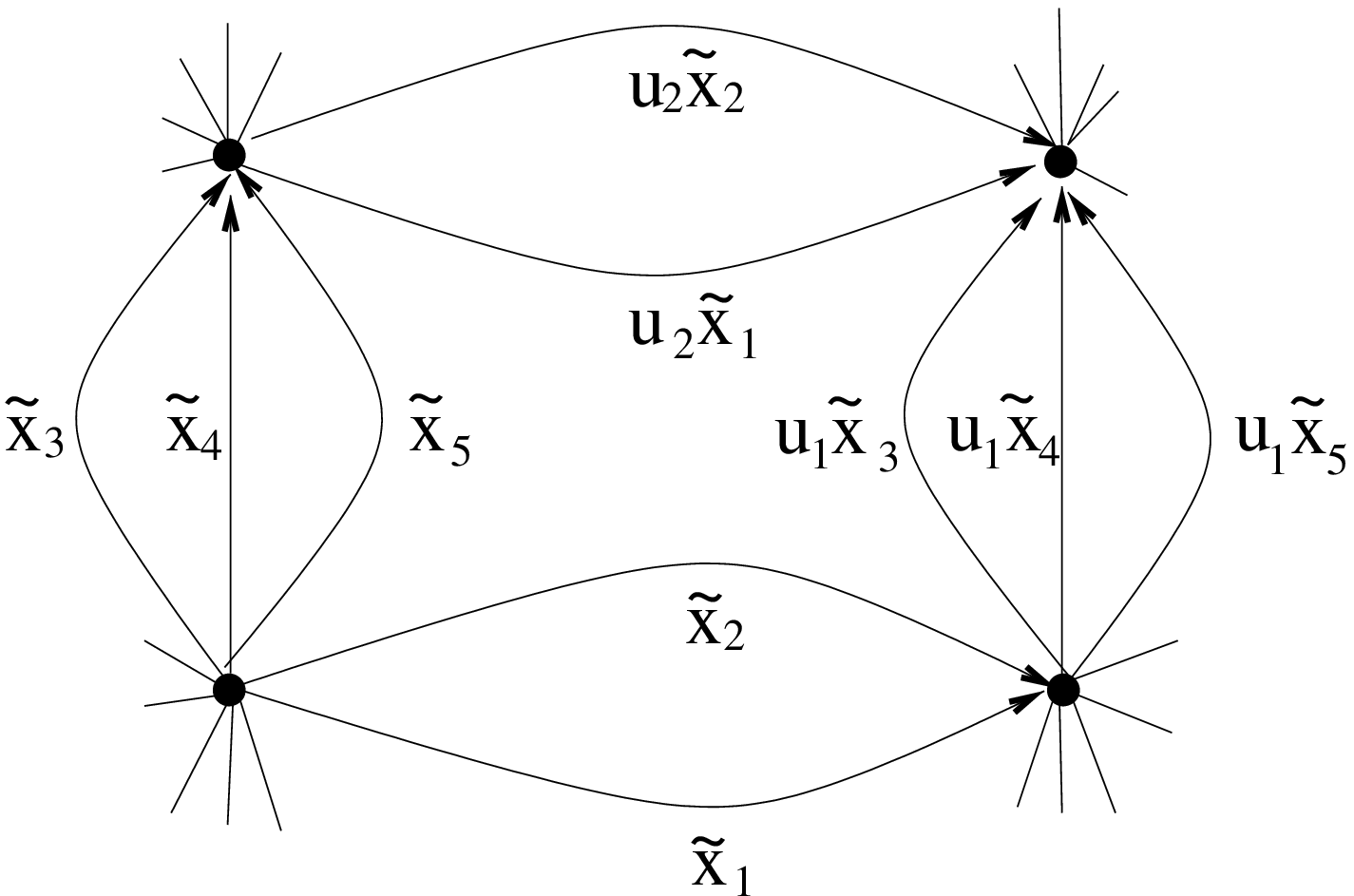}}
\bs

\centerline{{\bf Figure 3:} Portion of $\tilde P$}
\bs

We consider the subgroup $\L$ of $\Z^2$ generated by $u^3_1$ and 
$u^2_2$. The $1$-skeleton of $Q$ (which is the same as the $1$-skeleton of $\tilde P/\L$) is shown in
Figure 4.  Since $\L$ has generators parallel to $u_1, u_2$, it is easy to visualize the additional
$2$-cells that must be attached to $\tilde P/\L$ in order to build $Q$. In more general examples, the
boundaries of the new cells might wind around the graph several times.\bs 

\epsfxsize=3.4truein
\centerline{\epsfbox{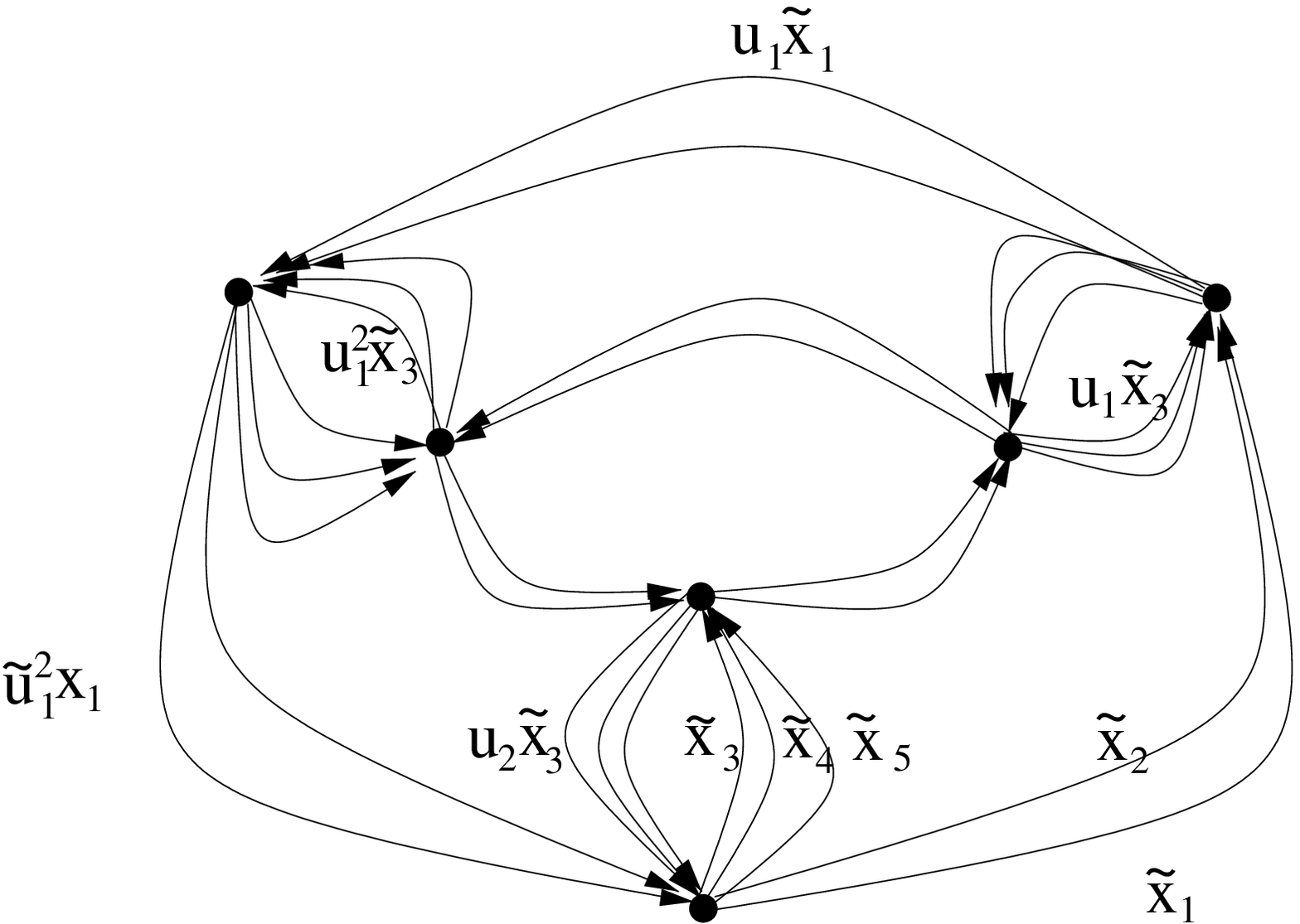}}
\bs
\centerline{{\bf Figure 4:} $1$-Skeleton of $Q$}
\bs

An element $x \in {\rm Fix}_\L(\s^B)$ can be represented by  $3 \times 2$ matrices $\b,\g, \d, \epsilon$ 
assigned to arcs corresponding to $x_2, x_3, x_4, x_5$, respectively. These matrices have entries in $\T$,
and are the restrictions of the elements of $\T^{\Z^d}$ to a fundamental region of $\Z^d/\L$. The element
$x$ is in ${\rm SFix}_\L(\s^B)$ if the column sums of $\b$ and the row sums of $\g, \d, \epsilon$ are all
zero.

An element of ${\rm Fix}_\L(\s')$ assigns additional  $3\times 2$ matrices $\zeta^\b, \zeta^\g, \zeta^d,
\zeta^\epsilon$. We can prescribe the coordinates $\zeta^\b_{0,0}, \zeta^\b_{0,1}$
arbitrarily; the other coordinates of $\zeta^\b$ are uniquely determined by these and 
$\b$. Similarly, the coordinates $\zeta^\g_{0,0}, \zeta^\g_{1,0}, \zeta^\g_{2,0}, \zeta^\d_{0,0},
\zeta^\d_{1,0}, \zeta^\d_{2,0},\zeta^\epsilon_{0,0},
\zeta^\epsilon_{1,0}, \zeta^\epsilon_{2,0}$ are arbitrary. We have ${\rm Fix}_\L(\s') \cong {\rm SFix}_\L(\s)
\oplus
\T^{11}$. 

We remark that for this link the Alexander polynomial is $(1-u_1)(1-u_2)$. Since the Mahler measure
of the polynomial is $1$, the orders $b_\L$ have zero exponential growth rate. \bs

\ni {\bf Proof of Corollary 2.2.} Since $\D_1/\D_2$
annihilates the Alexander module of any knot [{\bf Cr1}], it follows that 
$\s^N={\rm id}.$ From this we have ${\rm Fix}(\s^{r+N}) = {\rm
Fix}(\s^r)$. Recall that for any knot, ${\rm Fix}(\s^r)$ is
isomorphic to the dual of 
$H_1(M_r)$. Hence
$H_1(M_{r+N})\cong H_1(M_r)$, for every $r\ge 1$.\qed\bs

\ni {\bf 5. Coloring with nonabelian groups.} The coloring rule (3.1) generalizes in a natural way,
allowing one to replace $\T$ with an arbitrary topological group $\Sigma$. \bs

\ni {\bf Definition 5.1.} Assume that $D$ is a diagram of an oriented link $l=l_1 \cup \cdots l_d$
of $d$ components. A $\Sigma^{\Z^d}$-{\it coloring} of $D$ is an assignment of elements ({\it colors})
$\a, \b, \ldots \in \Sigma^{\Z^d}$ to the arcs of $D$ such that the condition 
$$\a\cdot \s_t\b = \g \cdot \s_{t'}\a \eqno(5.1)$$
is satisfied at any crossing. The colors $\a,\b,\g$ correspond to arcs that are described as in Definition
3.1. 

As before, if $D$ consists of $N$ arcs, then the set ${\rm Col}_{\Sigma, \Z^d}(D)$ of all
$\Sigma^{\Z^d}$-colorings of $D$ is a closed subspace of $[\Sigma^N]^{\Z^d}$ that is invariant under $\s_\m$
for each $\m \in {\Z^d}$. That ${\rm Col}_{\Sigma, \Z^d}(D)$ does not depend on the choice of diagram 
for $l$ follows immediately from the following. Let $\tilde E$ denote the maximal abelian cover
of the link exterior with projection $p: \tilde E \to E$, and let $*$ be a point of $E$. Let $\tilde *$
denote a fixed lift of $*$. Any covering automorphism of $\tilde E$ induces a
homeomorphism of the quotient space $\tilde E/p^{-1}(*)$, and hence induces an automorphism of 
$\pi_1(\tilde E/p^{-1}(*), \tilde *)$. By considering the adjoint action we obtain a homeomorphism
of the representation space ${\rm Hom}[\pi_1(E/p^{-1}(*), \tilde *), \Sigma]$. In this way we obtain a
$\Z^d$-action
$\s'$ on ${\rm Hom}[\pi_1(E/p^{-1}(*), \tilde *), \Sigma]$. A $\Z^d$-action on a topological space is
defined just as for $\Z^d$-action on a topological group, eliminating the requirement of a group structure
on the space; two $\Z^d$-actions, $\s$ acting on $X$ and $\s'$ acting on $X'$, are {\it
topologically conjugate} if there is a homeomorphism $\phi: X \to X'$ such that $\phi \circ \s_\m = 
\s'_\m \circ \phi,$ for each $\m \in \Z^d$.

\bs

\ni {\bf Proposition 5.2.} The $\Z^d$-actions on ${\rm Hom}[\pi_1(\tilde E/p^{-1}(*), \tilde *), \Sigma]$
and ${\rm Col}_{\Sigma, \Z^d}(D)$ are topologically conjugate. \bs

\ni {\bf Proof.} The quotient space $\tilde E/p^{-1}(*)$ has the same fundamental group as the quotient
complex $\tilde P/\tilde P^0$, where $\tilde P$ is defined in Section 4 and $\tilde P^0$
is the $1$-skeleton. There is a group
presentation for $\pi_1(\tilde P/\tilde P^0)$ in which the generators correspond to the edges of $\tilde P$;
lifts in $\tilde P$ of closed paths representing Wirtinger relators are closed paths in $\tilde P/\tilde
P^0$ representing the relators. Assignments of colors to the arcs of $D$, or equivalently to the
edges of $P$, such that the condition (5.1) holds at each crossing then correspond to homomorphisms from
$\pi_1(\tilde P/\tilde P^0)$ to
$\Sigma$. The correspondance defines a homeomorphism $\phi: {\rm Col}_{\Sigma, \Z^d}(D) \to {\rm
Hom}[\pi_1(\tilde E/p^{-1}(*), \tilde *),
\Sigma]$ such that $\phi\circ \s_\m = \s'_\m \circ \phi$, for each $\m \in \Z^d$.  \qed \bs

In view of this proposition we call ${\rm Col}_{\Sigma, \Z^d}(D)$ the {\it color $\Sigma^d$-shift}
of the link $l$, and we denote it by ${\rm Col}_{\Sigma, \Z^d}(l)$. 
 
The abelianization of $\pi_1(\tilde E/p^{-1}(*))$ is $H_1(\tilde P, \tilde P^0)$. It is
isomorphic to $H_1(\tilde E, p^{-1}(*))$, the Alexander module of the link. Hence we
propose that $\pi_1(\tilde E/p^{-1}(*))$ be called the {\it Alexander group} of the link; we
denote the group by ${\cal A}_l$. (Shortly after completing the
first draft of this paper, the authors discovered that the
Alexander group is a special case of the {\it derived group of  a
permutation representation}, introduced by R. Crowell [{\bf
Cr2}].)

It is much 
easier to write a presentation for the Alexander group than for the commutator subgroup of $\pi_1(S^3 -l)$.
One begins with families of generators $a_\m, b_\m,  c_\m,\dots$ $(\m \in \Z^d)$ corresponding to the arcs
of the diagram. Each crossing gives rise to a family of relations: a crossing such
as in Figure 1 imposes the relation $a_\m \cdot b_{\m+u_t} = c_\m \cdot a_{\m+ u_t}$. In the case of a knot,
when $d=1$, presentations of this sort are well known; J.~C. Hausmann and M. Kervaire [{\bf
HaKe}] termed them $\Z$-{\it dynamic}. A presentation for the Alexander group of a link such as the one we
have described might be called $\Z^d$-{\it dynamic}. 

The next proposition describes the relationship between the Alexander group ${\cal A}_l$ of a link $l$ and the 
commutator subgroup $G_l'$. We use the terminology of Section 4. Recall that $T$ is a maximal tree in the $1$-skeleton of
the cover $\tilde P$. \bs

\ni {\bf Proposition 5.3.} Let $l$ be an oriented link of $d$ components. The generators of ${\cal A}_l$
corresponding to the edges of $T$ freely generate a subgroup $F(E_T)$ of ${\cal A}_l$. Moreover, 
${\cal A}_l$ is isomorphic to the free product $G_l'*F(E_T)$. \bs

\ni {\bf Proof.} Let $C\tilde P^0$ denote the cone on the $0$-skeleton of $\tilde P$. The fundamental
group of $X= \tilde P \cup_{\tilde P^0} C\tilde P^0$ is isomorphic to ${\cal A}_l$. We can regard
$X$ as the union of $\tilde P$ and $T\cup _{\tilde P^0} C\tilde P^0$, which have contractible
intersection
$T$. An application of the Seifert van-Kampen theorem completes the argument. \qed \bs

\ni {\bf Example 5.4.} We examine two examples. Both are simple, but they highlight some of the 
advantages of working with the Alexander group rather the commutator subgroup of a link. 

\quad (i) Let $l$ be the trivial 2-component link. The group $G_l$ is free
on two generators. Choosing a diagram without crossings, we find that the Alexander group ${\cal A}_l$ is free
on  generators $a_{i,j}, b_{i,j}$, where $i, j$ range over $\Z$. The commutator subgroup $G'$ is also free,
but it does not admit a natural $\Z^2$-action by automorphisms as does ${\cal A}_l$.
\ss
\quad (ii) Next consider the link $l=2^2_1$, a Hopf link. The group $G_l$ is free abelian of rank 2. The
Alexander group ${\cal A}_l$ has presentation $\<a_{i,j}, b_{i,j}\mid a_{i,j}b_{i+1,j} = b_{i,j}a_{i, j+1}\>$.
In this example the
commutator subgroup $G_l'$ is trivial. \bs

\ni {\bf 6. Conclusion.} A possible direction for further inquiry involves links with zero Alexander
polynomial. 

The homology growth rate in Theorem 2.1 was computed
as the topological entropy of a $\Z^d$-action $\s'$. When the Alexander polynomial
of the link is zero, the entropy of $\s'$ can be shown to be infinite; in such a case we obtained no
information. However, we offer the following \ms

\ni {\bf Conjecture 6.1.} If $l$ is an oriented link of $d$ components, then 
$$\lim_{\<\L\>\to \infty}{1 \over {|\Z^d/\L|}}\log|TH_1(M_\L(l))|= \log {\bf M}(\D_i),$$
where $\D_i$ is the first {\sl nonzero} Alexander polynomial of the link. \ms


\centerline{\bf References.} 
\bs
\baselineskip=12 pt
\item{[{\bf Ah}]} L.~V. Ahlfors, Complex analysis, 2nd edition, McGraw-Hill, New York, 1966.
\ss
\item{[{\bf Ba}]} A.~Baker, ``The theory of linear forms in logarithms,'' in Transcendence Theory,
Advances and Applications, Academic Press, 1977. 
\ss
\item{[{\bf Cr1}]} R.~H.Crowell, ``The annihilator of a knot module,'' {\sl Proc.\ Amer.\
Math.\  Soc.\ \bf 15} (1964), 696 -- 700.
\ss
\item{[{\bf Cr2}]} R.~H.Crowell, ``The derived group of a permutation
representation,'' {\sl Advances\ in\ Math.\ \bf 53}(1984), 99 -- 124.
\ss
\item{[{\bf EvWa}]} G. Everest, T. Ward, Heights of polynomials and entropy in algebraic dynamics,
Springer-Verlag, London, 1999. 
\ss
\item{[{\bf Fo}]} R.~H. Fox, ``A quick trip through knot theory,'' in Topology of $3$-Manifolds and
Related Topics (edited by M.~K. Fort), Prentice-Hall, NJ (1961), 120--167. 
\ss
\item{[{\bf Ge}]} A.~O. Gel'fond, ``On the approximation of transcendental numbers  by algebraic
integers,'' {\sl Dokl.\ Akad.\ Navk.\ SSSR\ \bf 2} (1935), 177 -- 182.
\ss
\item{[{\bf Go}]} C.~McA. Gordon, ``Knots whose branched coverings have periodic homology,''
{\sl Trans.\ Amer.\ Math.\ Soc.\ \bf 168} (1972), 357 -- 370.
\ss
\item{[{\bf GoSh}]} F. Gonz\'alez-Acu\~na and H. Short, ``Cyclic branched coverings of knots and 
homology spheres,'' {\sl Revista\ Math.\ \bf 4} (1991), 97 -- 120.
\ss
\item{[{\bf HaKe}]} J.~C. Hausmann and M. Kervaire, ``Sous-groupes d\'eriv\'es des groupes
des noeuds,'' {\sl L'Enseign.\ Math.\ \bf 24} (1978), 111 -- 123.
\ss
\item{[{\bf Le}]} D.~H. Lehmer, ``Factorization of certain cyclotomic functions,'' {\sl Annals\ of\
Math.\ \bf 34} (1933), 461 -- 479. 
\ss
\item{[{\bf Li}]} W.~B. Lickorish, An Introduction to Knot Theory, Springer-Verlag, Berlin, 1997. 
\ss
\item{[{\bf LiScWa}]} D. Lind, K. Schmidt and T. Ward, ``Mahler measure and entropy for commuting
automorphisms of compact groups,'' {\sl Invent.\ Math.\ \bf 101} (1990), 593--629.
\ss
\item{[{\bf Ma1}]} K. Mahler, ``An application of Jensen's formula to polynomials,'' {\sl Mathematika\ 
\bf 7} (1960), 98 -- 100.
\ss
\item{[{\bf Ma2}]} K. Mahler, ``On some inequalities for polynomials in several variables,'' 
{\sl J.\ London\ Math.\ Soc.\ \bf 37}(1962), 341 -- 344. 
\ss
\item{[{\bf Ri}]} R. Riley, ``Growth of order of homology of cyclic branched covers of knots,''
{\sl Bull.\ London\ Math.\ Soc.\ \bf 22} (1990), 287 -- 297.
\ss
\item{[{\bf Ro}]} D. Rolfsen, Knots and Links, Publish or Perish, Berkeley, CA, 1976. 
\ss
\item{[{\bf Sa}]} M. Sakuma, ``Homology of abelian coverings of links and spatial graphs,''
{\sl Canad.\ J.\ Math.\ \bf 47} (1995), 201 -- 224.
\ss
\item{[{\bf Sc}]} K. Schmidt, Dynamical Systems of Algebraic Origin, Birkh\"auser Verlag, Basel, 1995.
\ss
\item{[{\bf SiWi1}]} D.~S. Silver and S.~G. Williams, ``Generalized $n$-colorings of links,'' in 
Knot Theory, Banach Center Publications 42, Warsaw 1998, 381--394. 
\ss
\item{[{\bf SiWi2}]} D.~S. Silver and S.~G. Williams, ``Coloring link diagrams with a continuous
palette,'' {\sl Topology\ \bf 39} (2000), 1225 -- 1237.
\ss
\item{[{\bf SiWi3}]} D.~S. Silver and S.~G. Williams, ``Knot invariants from symbolic dynamical
systems,'' {\sl Trans.\ Amer.\ Math.\ Soc.\ \bf 351} (1999), 3243--3265.
\ss
\item{[{\bf SiWi4}]} D.~S. Silver and S.~G. Williams, ``Mahler measure of Alexander
polynomials,'' preprint, 2000. 
\ss
\item{[{\bf Sp}]} E.~H. Spanier, Algebraic Topology, McGraw-Hill Book Co., New York, 1966.
\ss
\bs
\item{} Dept. of Mathematics and Statistics, Univ. of South Alabama, Mobile, AL  36688-0002
e-mail: silver@mathstat.usouthal.edu, williams@mathstat.usouthal.edu

\end